\documentclass[12pt]{article}
\usepackage{amsmath}
\usepackage{amsthm}
\usepackage{amssymb}
\usepackage{eucal}
\usepackage{graphicx,graphics}
\numberwithin{equation}{section}
\usepackage{graphicx,graphics}
\usepackage{pdfsync}
\def \dis {\displaystyle}

\def \into {\int_\Omega}
\def \confai {-\kern -.5em\rightharpoonup}
\def \cqfd {\hfill$\Box$}

\def \ep {\varepsilon}

\def \Om {\Omega}

\def \NN {\mathbb N}

\def \RR {\mathbb R}
\def \beq {\begin{equation}}
\def \eeq {\end{equation}}
\def \ba {\begin{array}}
\def \ea {\end{array}}
\def \bs {\bigskip}

\def \ecart {\noalign{\medskip}}
\newtheorem{Thm}{Theorem}[section]
\newtheorem{Cor}[Thm]{Corollary}

\newtheorem{Adef}[Thm]{Definition}

\newtheorem{Arem}[Thm]{Remark}
\newenvironment{Rem}{\begin{Arem}\rm}{\end{Arem}}
\newtheorem{Aexa}[Thm]{Example}

\newtheorem{Anot}[Thm]{Notation}

\def \refe #1.{(\ref{#1})}
\def \reff #1.{figure~\ref{#1}}
\def \refs #1.{Section~\ref{#1}}
\def \refss #1.{Subsection~\ref{#1}}
\def \refD #1.{Definition~\ref{#1}}
\def \refT #1.{Theorem~\ref{#1}}
\def \refL #1.{Lemma~\ref{#1}}
\def \refC #1.{Corollary~\ref{#1}}
\def \refP #1.{Proposition~\ref{#1}}
\def \refPt #1.{Properties~\ref{#1}}
\def \refR #1.{Remark~\ref{#1}}
\def \refE #1.{Example~\ref{#1}}
\def \refN #1.{Notation~\ref{#1}}
\title{An elliptic problem in dimension $N$ with a varying drift term bounded in $L^N$.}
\author{
\footnotesize
\centerline{\begin{tabular}{cccc}
 \normalsize Juan Casado-D\'{\i}az
\\
 Dpto. de Ecuaciones Diferenciales y An\'alisis Num\'erico
\\
 Universidad de Sevilla
\\
jcasadod@us.es
\end{tabular}}
}
\begin{document}
\maketitle
\begin{abstract}
The present paper is devoted to study the asymptotic behavior of  a sequence of linear elliptic equations with a varying drift term, whose coefficients are just bounded in $L^N(\Om)$, with $N$ the dimension  of the space. It is known that there exists a unique solution for each of these problems in the Sobolev space $H^1_0(\Om)$. However, because the operators are not coercive, there is no uniform estimate of the solutions  in this space. We use some estimates in \cite{boc1}, and a regularization obtained by adding a small nonlinear first order term, to pass to the limit in these problems.\end{abstract}
\par\bs\noindent
\par\bs\noindent
\section{Introduction.}
For a bounded open set  $\Om\subset\RR^N$, we are interested  in passing to the limit in a sequence on elliptic equations with a varying  drift term, whose coefficients are just bounded in $L^N(\Om)^N$, ($L^p(\Om)^N$, $p>2$ if $N=2$). The problem is written as
\beq\label{pbint1} \left\{\ba{l}\dis -{\rm div}(A\nabla u_n+E_nu_n)=f_n\ \hbox{ in }\Om\\ \ecart\dis u_n\in H^1_0(\Om),\ea\right.\eeq
with $A\in L^\infty(\Om)^{N\times N}$ satisfying the usual uniform ellipticity condition, and $f_n$ bounded in $H^{-1}(\Om)$. Thanks to Sobolev's inequality,
the integrability assumption on $E_n$ is the weaker one to get the  first order term well defined in $H^{-1}(\Om)$.
It is known that for every $n\in\NN$ problem (\ref{pbint1}) has a unique solution (\cite{boc2}, \cite{bo-icd}, \cite{GiTr}), however, the problem is known to be not coercive and thus, there is no estimate for $u_n$ in $H^1_0(\Om)$.  We recall some of these results in Section \ref{SeEx}.\par
L. Boccardo found in \cite{boc1} an estimate for 
$\ln (1+|u_n|)$ in $H^1_0(\Om)$ depending only on $\|f_n\|_{H^{-1}(\Om)}$ and $\|E_n\|_{L^2(\Om)^N}$.  As a consequence of this result we get that the measure of the sets $\{|u_n|>k\}$ tends to zero when $k$ tends to infinity uniformly in $n$. Thanks to this result L. Boccardo proved in \cite{boc3} that $f_n$ converging weakly to some $f$ in $H^{-1}(\Om)$, $E_n$ converging weakly to $E_0$ in $L^N(\Om)$ and $|E_n|^N$  equi-integrable, imply that the solution $u_n$ of (\ref{pbint1}) converges weakly in $H^1_0(\Om)$ to the solution $u_0$ of
$$ \left\{\ba{l}\dis -{\rm div}(A\nabla u_0+E_0u_0)=f\ \hbox{ in }\Om\\ \ecart\dis u_0\in H^1_0(\Om).\ea\right.$$
In the case where $|E_n|^N$ is not equi-integrable, we do not have an estimate for $u_n$ in $H^1_0(\Om)$ and then the proof is more involved. The ideas in \cite{boc1}  imply (see \cite{boc3}) that $\ln^q (1+|u_n|)$ is bounded in $H^1_0(\Om)$ for any $q>1$. Using also the approximation of $u_n$ given by the solution of a perturbation of (\ref{pbint1}) by a nonlinear zero order term (see (\ref{defwnd}) below), we still manage to pass to the limit in (\ref{pbint1}) but, instead of the weak convergence in $H^1_0(\Om)$, we only get
$$\ln^q(1+|u_n|)\,{\rm sgn}(u_n)\rightharpoonup \ln^q(1+|u_0|)\,{\rm sgn}(u_0)\ \hbox{ in }H^1_0(\Om),\quad\forall\, q\geq 1.$$
This is the main result of the paper, which is proved in Theorem \ref{Thgho}.
\par 
Another problem related to (\ref{pbint1}) is given by its adjoint formulation
\beq\label{ptint2}\left\{\ba{l}\dis -{\rm div}(A^T\nabla u_n)+E_n\cdot \nabla u_n=f_n\ \hbox{ in }\Om\\ \ecart\dis u_n\in H^1_0(\Om),\ea\right.\eeq
Some results about the asymptotic behavior of the solutions of this problem have been obtained in \cite{boc3}.
 Namely, assuming the equi-integrability condition on $|E_n|$ and reasoning by duality, we deduce from the results stated above for problem (\ref{pbint1})  that the solutions of (\ref{ptint2}) are bounded in $H^1_0(\Om)$. Moreover, the results in  \cite{BoMu} prove that $u_n$ is compact in $W^{1,q}_0(\Om)$ for $1\leq q<2$. Then, it is immediate to pass to the limit in this problem. In the case where $|E_n|^N$ is not equi-integrable,  we do not have any estimate for $u_n$ and then we are not able to passing to the limit in (\ref{ptint2}). However, adding a sequence of zero order terms $a_nu_n$ with $a_n$ satisfying
$$a_n\geq \gamma\ \hbox{ a.e. in }\Om,$$
for some $\gamma>0$, and assuming $f_n$ bounded in $L^\infty(\Om)$, it is proved in \cite{boc3} that $u_n$ is bounded in $H^1_0(\Om)\cap L^\infty(\Om)$. This allows us to pass to the limit in (\ref{ptint2}). \par
The homogenization of a sequence of elliptic PDE's with a singular term of first order has also been carried out in other papers. In this sense we refer to \cite{BoCD} where it is considered problem (\ref{pbint1}) with the matrix function $A$ also depending on $n$, the sequence $E_n$  converging weakly in $L^2(\Om)^N$, and ${\rm div}\,E_n$ converging strongly in $H^{-1}(\Om)$. The definition of solution in this case is related to the definition of entropy or renormalized solution (see e.g. \cite{BBGGPV}, \cite{BoGa}, \cite{DMMOP}). In \cite{BrCa} and \cite{BrGe} it is considered the case of a first order term of the form
$$E_n\nabla u_n+{\rm div}(E_nu_n).$$
Since this term is skew-symmetric we can obtain an estimate in $H^1_0(\Om)$ which is independent of $E_n$ (we refer to \cite{BoCDOr} for a related existence result). Assuming $E_n$ just bounded in $L^2(\Om)^N$ it is proved that the limit problem contains a new term of zero order. This is related to the results
 obtained in \cite{Tar} for the Stokes equation with an oscillating Coriolis force. We also refer to \cite{BrCa2}, \cite{BrCa3} for related results in the case of the evolutive elastic system submitted to an oscillating magnetic field. Now, the limit problem is nonlocal in general.
\section{Some reminders about  elliptic problems with a drift or convection term.}\label{SeEx}
For a bounded open set $\Om\subset\RR^N$, $N\geq 2$, a matrix function $A\in L^\infty(\Om)^{N\times N}$, such that there exists $\alpha>0$ satisfying
\beq\label{elip}A(x)\xi\cdot\xi\geq \alpha|\xi|^2,\quad\forall\,\xi\in \RR^N,\ \hbox{ a.e. }x\in \Om,\eeq
two measurable functions $E:\Om\to\RR^N$, $a:\Om\to\RR$, $a\geq 0$ a.e. in $\Om$, and a distribution $f\in H^{-1}(\Omega)$,
we recall in this section, some results about the existence and uniqueness of solution for problems
\beq\label{PbEDP}\left\{\ba{l} {\cal A}u=f\ \hbox{ in }\Om\\ \ecart\dis
u\in H^1_0(\Omega),\ea\right.\qquad \left\{\ba{l} {\cal A}^\ast u=f\ \hbox{ in }\Om\\ \ecart\dis
u\in H^1_0(\Omega),\ea\right.\eeq
with
$${\cal A}u=-{\rm div}(A\nabla u+Eu)+au,\qquad {\cal A}^\ast u=-{\rm div}(A^T\nabla u)+E\cdot\nabla u+au,$$
Observe that due to  Sobolev imbedding theorem, in order to have the terms ${\rm div}(Eu)$, $E\cdot\nabla u$ and $au$ in $H^{-1}(\Omega)$, when $u$ is in $H^1_0(\Omega)$, we need to assume
\beq\label{regE} \left\{\ba{cc}E\in L^N(\Om)^N&\hbox{ if }N>2\\ \ecart\dis 
E\in L^p(\Om)^N,\ p>2 &\hbox{ if }N=2,\ea\right.\qquad \left\{\ba{cc}a\in L^{N\over 2}(\Om)&\hbox{ if }N>2\\ \ecart\dis 
a\in L^q(\Om)^N,\ q>1 &\hbox{ if }N=2,\ea\right. \eeq
In this case, ${\cal A}$ and ${\cal A}^\ast$ are continuous linear operators from $H^1_0(\Om)$ into $H^{-1}(\Om)$ and ${\cal A}^\ast$ is the adjoint operator of ${\cal A}$.\par
Assuming the further assumption 
\beq\label{hipmf}E\in L^p(\Om)^N,\quad p>N,\eeq
the compact imbedding of $H^1_0(\Om)$ into $L^{2p\over p-2}(\Omega)$ proves that for every $\ep>0$, there exists $C_\ep>0$ such that
\beq\label{concom}\|u\|_{L^{2p\over p-2}(\Omega)}\leq \alpha\ep\|\nabla u\|_{L^2(\Om)^N}+C_\ep\|u\|_{L^2(\Om)},\quad \forall\, u\in H^1_0(\Om).\eeq
Thus, H\"older's and Young's inequalities imply
$$\ba{l}\dis \langle {\cal A}u,u\rangle_{H^{-1}(\Om),H^1_0(\Om)}=\langle {\cal A}^\ast,u\rangle_{H^{-1}(\Om),H^1_0(\Om)}
=\into \big(A\nabla u\cdot\nabla u\,dx+uE\cdot \nabla u+au^2\big)dx\\ \ecart\dis\geq \alpha\|\nabla u\|_{L^2(\Om)^N}^2-\|E\|_{L^p(\Om)^N}\|\nabla u\|_{L^2(\Om)^N}\|u\|_{L^{2p\over p-2}(\Omega)}\\ \ecart\dis
\geq \alpha\Big({1\over 2}-\ep\|E\|_{L^p(\Om)^N}\Big)\|\nabla u\|^2_{L^2(\Om)^N}
-{C_\ep^2\|E\|_{L^p(\Om)^N}^2\over 2\alpha}\|u\|^2_{L^2(\Om)}.\ea$$
Therefore, taking $\ep\|E\|_{L^p(\Om)^N}<1/2$, we deduce from Lax-Milgram theorem that replacing $a$
by $a+\gamma$, with
$$\gamma\geq {C_\ep^2\|E\|_{L^p(\Om)^N}^2\over 2\alpha},$$
there exists a unique solution for both problems in (\ref{PbEDP}).
The compactness of  $({\cal A}+\gamma I)^{-1}$, considered as an operator in $L^2(\Om)$, allows then to use Fredholm theory to deduce that the existence and uniqueness of solutions for both problems in (\ref{PbEDP}) and every $f\in H^{-1}(\Om)$ is equivalent to the uniqueness of solutions for one of them. \par
In \cite{GiTr}, Theorem 8.1, it is proved that  problem
$$\left\{\ba{l} {\cal A}^\ast u=f\ \hbox{ in }\Om\\ \ecart\dis
u-\phi\in H^1_0(\Omega),\ea\right.$$
with $\phi\in H^1(\Om)$ satisfies the weak maximum principle so that the second problem in (\ref{PbEDP}) has at most one solution (in \cite{GiTr} it is assumed $E\in L^\infty(\Om)^N$, but it is immediate to check that the same proof works for $E$ just satisfying (\ref{regE})). \par
We observe however that although the above result implies the existence and continuity of the operators ${\cal A}^{-1}$ and $({\cal A}^\ast)^{-1}$, it does not provide any estimate for the norm of these operators and then on the solutions of both problems in (\ref{PbEDP}).\par
If we assume $N>2$ and $E\in L^N(\Omega)^N$, the above reasoning fails because (\ref{concom}) with $p=N$ does not hold in general. The existence of solutions in this case can be obtained from the
 following result due to L. Boccardo (\cite{boc1}, \cite{boc3}). We recall that the truncate function at height $k>0$ is defined as
\beq\label{deTk} T_k(s)=\left\{\ba{cl} -k &\hbox{ if }s<-k\\ \ecart\dis s &\hbox{ if }-k\leq s\leq k\\ \ecart\dis k &\hbox{ if }s>k.\ea\right.\eeq
\begin{Thm} \label{Solde} For every $f\in H^{-1}(\Om)$, $E\in L^2(\Om)^N$, and $a\in L^1(\Om)$, $a\geq 0$ a.e. in $\Om$, there exists an entropy solution of
\beq\label{pbD1}\left\{\ba{l}\dis {\cal A}u=f\ \hbox{ in }\Om\\ \ecart\dis u=0\ \hbox{ on }\partial\Om,\ea\right.\eeq
in the following sense
\beq\label{defres}\left\{\ba{l}\dis T_k(u)\in H^1_0(\Om),\ \forall\, k>0,\quad \log(1+|u|)\in H^1_0(\Om)\\ \ecart\dis
\into \Big((A\nabla u+Eu)\cdot\nabla T_k(u-\varphi)+a uT_k(u-\varphi)\Big)\,dx\leq \langle f,T_k(u-\varphi)\rangle_{H^{-1}(\Om),H^1_0(\Om)},\\ \ecart\dis \forall\, \varphi\in H^1_0(\Om)\cap L^\infty(\Om).\ea\right.\eeq
Moreover, there exists $C>0$, independent of $f$ and $E$ such that
\beq\label{espriL}\|\log(1+|u|)\|_{H^1_0(\Om)}\leq C\big(\|f\|_{H^{-1}(\Om)}+\|E\|_{L^2(\Om)^N}\big).\eeq
If there exists $\gamma>0$ such that
\beq\label{elipa}a\geq\gamma\ \hbox{ a.e. in }\Om,\eeq
and $f$ belongs to $L^1(\Om)$, such solution is also in $L^1(\Om)$ and satisfies
\beq\label{acoL1} \|u\|_{L^1(\Om)}\leq {1\over \gamma}\|f\|_{L^1(\Om)}.\eeq
\end{Thm}
Let us prove that the function $u$ given by the previous theorem is in fact a distributional solution of the first problem in (\ref{PbEDP}) when $E$ and $a$ satisfy (\ref{regE}). This is given by the following theorem
\begin{Thm} \label{thmexg} Assume that the functions $E$ and $a$ in Theorem \ref{Solde} satisfy (\ref{regE}), then the solution $u$ of (\ref{pbD1}) is in $H^1_0(\Om)$ and satisfies the elliptic equation in the distributional sense. Moreover, if there exists $\gamma>0$ such that (\ref{elipa}) holds, then this solution is unique.
\end{Thm}\par\noindent
{\bf Proof.} For $k,m>0$, we take  $\varphi=T_m(u)$ in (\ref{defres}). This gives
$$\int_{\{m<|u|<k+m\}}\hskip-15pt\big((A\nabla u+Eu)\cdot\nabla u\big)dx+  \into au\,T_k(u-T_m(u))dx
\leq \langle f,T_k(u-T_m(u))\rangle_{H^{-1}(\Om),H^1_0(\Om)}.$$
Using in this equality that
$$\left|\int_{\{m<|u|<k+m\}}\hskip-15pt Eu\cdot\nabla u\,dx\right|\leq \int_{\{m<|u|<k+m\}}\hskip-15pt |E|\big(m+|T_k(u-T_m(u))|\big)|\nabla u|dx,$$
combined with (\ref{elip}), Sobolev's imbedding theorem and H\"older's inequality we deduce the existence of $C_S>0$, depending only on $N$ (and $|\Om|$ if $N=2$), such that
\beq\label{acotH1}\ba{l}\dis \|T_k(u-T_m(u))\|^2_{H^1_0(\Om)}\big(\alpha-C_S\|E\|_{L^r(\{m<|u|\})^N}\big)\\ \ecart\dis \leq
\big(\|f\|_{H^{-1}(\Om)}+mC_S\|E\|_{L^r(\Om)^N}\big)\|T_k(u-T_m(u))\|_{H^1_0(\Om)},\ea\eeq
with $r=p$ if $N=2$, $r=N$ if $N>2$. \par
By (\ref{espriL}), we can take $m$ sufficiently large to get
$$C_S\|E\|_{L^r(\{m<|u|\})}^N<\alpha.$$
Letting $m$ fixed satisfying this inequality and taking $k$ tending to infinity, we deduce that 
$u-T_m(u)$ belongs to $H^1_0(\Om)$. Since $T_m(u)$ also belongs to $H^1_0(\Om)$, then $u$ belongs to $H^1_0(\Om)$.\par
Once we know that $u$ belongs to $H^1_0(\Om)$, we can take $k$ tending to infinity in (\ref{defres}) to deduce
$$\into \Big((A\nabla u+Eu)\cdot\nabla (u-\varphi)+a u(u-\varphi)\Big)\,dx\leq \langle f,u-\varphi\rangle_{H^{-1}(\Om),H^1_0(\Om)},$$
for every $\varphi\in H^1_0(\Om)\cap L^\infty(\Om).$ Since $H^1_0(\Om)\cap L^\infty(\Om)$ is dense in $H^1_0(\Om)$, the result holds for $\varphi$ just in $H^1_0(\Om)$. Replacing $\varphi$ by $u+\varphi$
we deduce that $u$ is a distributional solution of  ${\cal A}u=f$ in $\Omega$.\par
The fact that the renormalized solutions are distributional solutions for $E$ satisfying (\ref{regE}) allows us to reason by linearity to deduce that they are unique if and only if the unique entropy solution of (\ref{pbD1}) with $f=0$ is the null function. When (\ref{elipa}) holds this follows from (\ref{acoL1}). \cqfd\par\medskip
Theorem \ref{Solde} proves that the first problem in (\ref{PbEDP}) always has at least a solution 
and it is unique if $a$ satisfies (\ref{elipa}). As above, this allows us to use Fredholm theory to deduce
\begin{Cor} \label{Tges} Assume that $E$ satisfies (\ref{regE}), then for every $f\in H^{-1}(\Om)$  both problems in (\ref{PbEDP}) have a unique solution. Moreover this solution depends continuously on $f$.\end{Cor}
\section{Passing to the limit with a varying drift term.} \label{sec3}
In this section, for a bounded open set $\Om\subset\RR^N,$ $N\geq 2$, a sequence of vector measurable functions 
 $E_n:\Om\to\RR^N$  such that there exists  $E_0:\Om\to\RR^N$ satisfying
\beq\label{assEn} \left\{\ba{cc}E_n\rightharpoonup E_0\ \hbox{ in }L^N(\Om)^N&\hbox{ if }N>2\\ \ecart\dis 
E_n\rightharpoonup E_0\ \hbox{ in }L^p(\Om)^N,\hbox{ for some }p>2 &\hbox{ if }N=2,\ea\right.\eeq 
we are interested  in the asymptotic behavior of the solutions of 
\beq\label{pbhom} \left\{\ba{l}\dis -{\rm div}(A\nabla u_n+E_nu_n)=f_n\ \hbox{ in }\Om\\ \ecart\dis u_n\in H^1_0(\Om),\ea\right.\eeq
where $A\in L^\infty(\Om)^{N\times N}$ satisfies (\ref{elip}) and $f_n\in H^{-1}(\Om)$ is such that there exists $f\in H^{-1}(\Om)$ satisfying
\beq\label{confn} f_n\rightharpoonup f\ \hbox{ in }H^{-1}(\Om).\eeq
As we recalled in the previous section, problem (\ref{pbhom}) has a unique solution but we do not know if it is bounded in $H^1_0(\Om)$. Using the estimates for the solutions of this problem given in \cite{boc1} let us prove
\begin{Thm} \label{Thgho} Assume $\Om$ a bounded open set of $\RR^N$, $N\geq 2$, and $E_n$ a sequence of vector measurable functions in $\Om$\ such that there exists $E_0$ satisfying (\ref{assEn}).
Then, for every sequence $f_n$ which satisfies (\ref{confn}), the sequence of solutions $u_n$ of (\ref{pbhom})
satisfies
\beq\label{codTm} T_m(u_n)\rightharpoonup T_m(u_0)\  \hbox{ in }H^1_0(\Om),\ \forall\, m>0,\eeq
\beq\label{codln} \ln^q(1+|u_n|)\rightharpoonup \ln^q(1+|u_0|)\ \hbox{ in }H^1_0(\Om),\ \forall\, q>0,\eeq
with $u$ the unique solution of
\beq\label{pbli} \left\{\ba{l}\dis -{\rm div}(A\nabla u_0+E_0u_0)=f\ \hbox{ in }\Om\\ \ecart\dis u_0\in H^1_0(\Om).\ea\right.\eeq
Moreover, if one of the following assumptions hold:
\begin{itemize}
\item $N=2$.
\item $N>2$, $|E_n|^N$ is equi-integrable,
\end{itemize}
then the convergence holds in the weak topology of $H^1_0(\Om)$.
\end{Thm}
\par\noindent
{\bf Proof.} \par\medskip\noindent
{\it Step 1.} We start getting some estimates for the solutions of (\ref{pbhom}). They are based on \cite{boc1}.
\par
For $t\geq 2$, we use as test function in (\ref{pbhom}) $v=\phi(u_n)$, with $\phi$ defined by
$$\phi(s)=\int_0^s{\ln^t(1+ |\rho|)\over (1+|\rho|)^2}d\rho,\quad \forall\,s\in\RR.$$
We get
$$\into \big(A\nabla u_n+E_n u_n)\cdot\nabla u_n{\ln^t(1+ |u_n|)\over (1+|u_n|)^2}dx=\langle f_n,\phi(u_n)\rangle.$$
Using  (\ref{elip}), Young's and H\"older's inequalities, and $\phi'\in L^\infty(\RR)$, we deduce the existence of $C>0$ depending on $\alpha, t, N$ ($\alpha$, $t,$ $p$, $|\Om|$ if $N=2$) such that
\beq\label{Est1}\ba{l}\dis \into \Big|\nabla \Big(\ln^{t+2\over 2}(1+|u_n|)\,{\rm sgn}(u_n)\Big)\Big|^2dx
={(t+2)^2\over 4}\into |\nabla u_n|^2{\ln^t(1+ |u_n|)\over (1+|u_n|)^2}dx\\ \ecart\dis\leq C\into |E_n|^2|u_n|^2{\ln^t(1+ |u_n|)\over (1+|u_n|)^2}dx+C\|f_n\|_{H^{-1}(\Om)}^2\\ \ecart\dis\leq C\into |E_n|^2\ln^t(1+ |u_n|)dx+C\|f_n\|_{H^{-1}(\Om)}^2\\ \ecart\dis
\leq C\left(\into |E_n|^rdx\right)^{2\over r}\left(\into  \ln^{tr\over r-2}(1+ |u_n|)dx\right)^{r-2\over r}+C\|f_n\|_{H^{-1}(\Om)}^2,\ea\eeq
with 
\beq\label{defir} r=\left\{\ba{cl} p &\hbox{ if }N=2\\ \ecart\dis r=N &\hbox{ if }N>2.\ea\right.\eeq
\par 
By Sobolev's inequality, we also have
$$ \left(\into \ln^{(t+2)\mu\over 2}(1+|u_n|)dx\right)^{2\over \mu}\leq C_S\into \Big|\nabla \Big(\ln^{t+2\over 2}(1+|u_n|)\,{\rm sgn}(u_n)\Big)\Big|^2dx,$$
with $1\leq \mu$, if $N=2$ and $1\leq\mu\leq 2N/(N-2)$, if $N>2$. Choosing
$$1<\mu:={2tr\over (t+2)(r-2)}<{2r\over r-2}$$
and using Young's inequality, we deduce (for another constant $C>0$), after simplifying equal terms,
$$\left(\into  \ln^{tr\over r-2}(1+ |u|)dx\right)^{(t+2)(r-2)\over tr}\leq C\left(\into |E_n|^rdx\right)^{t+2\over r}+C\|f_n\|^2_{H^{-1}(\Om)}.$$
Replacing this estimate in the right-hand side of (\ref{Est1}), we finally get
\beq\label{estlog}\into \Big|\nabla \Big(\ln^{t+2\over 2}(1+|u_n|)\,{\rm sgn}(u_n)\Big)\Big|^2dx\leq 
C\left(\into |E_n|^rdx\right)^{t+2\over r}+C\|f_n\|^2_{H^{-1}(\Om)}.\eeq
Since $t$ can be taken as large as we want, we conclude 
\beq\label{estlog2} \ln^q(1+|u_n|)\,{\rm sgn}(u_n)\ \hbox{ bounded in }H^1_0(\Om),\quad\forall\,q>0.\eeq
In particular, $T_m(u_n)$ is bounded in $H^1_0(\Om)$, for every $m\in\NN$. These estimates and the Rellich-Kondrachov's compactness theorem prove the existence of a subsequence of $n$, still denoted by $n$, and a measurable function $u$ such that for every $m,q>0$, we have
\beq\label{codTm0} T_m(u_n)\rightharpoonup T_m(u)\hbox{ in }H^1_0(\Om),\quad \ln^q(1+|u_n|)\rightharpoonup \ln^q(1+|u|)\hbox{ in }H^1_0(\Om).\eeq
\par\medskip
Now, we have to prove that $u=u_0$.\par\medskip\noindent
{\it Step 2.} Let us first consider the case $N=2$ or $N\geq 3$, $|E_n|^N$ equi-integrable. It has been first carried out in \cite{boc3}. Indeed, since for  $N=2$, $p$ in (\ref{assEn}) is any number bigger than 2, the problem reduces to assume $|E_n|^r$ equi-integrable with $r$ given by (\ref{defir}). Taking as test function in (\ref{pbhom}) $u_n-T_m(u_n)$, we deduce
$$\int_{\{m<|u_n|\}}\big(A\nabla u_n+E_nu_n\big)\cdot\nabla u_n\,dx=\langle f_n,u_n-T_m(u_n)\rangle_{H^{-1}(\Om), H^1_0(\Om)},$$
which similarly to (\ref{acotH1}) implies
$$ \ba{l}\dis \|u_n-T_m(u_n)\|^2_{H^1_0(\Om)}\big(\alpha-C_S\|E_n\|_{L^r(\{m<|u_n|\})^N}\big)\\ \ecart\dis \leq
\big(\|f_n\|_{H^{-1}(\Om)}+m C_S\|E_n\|_{L^2(\Om)^N}\big)\|u_n-T_m(u_n)\|_{H^1_0(\Om)}.\ea$$
Thanks to (\ref{estlog2}), we have
$$\lim_{m\to\infty}\sup_{n\in\NN} \big|\{|u_n|>m\big\}\big|=0,$$
which combined with the equi-integrability of $|E_n|^r$ implies the existence of $m\in\NN$ such that $\|E_n\|_{L^r(\{m<|u_n|\})^N}<\alpha/(2C_S)$, for every $n\in\NN$. Therefore $u_n-T_m(u_n)$ is bounded in $H^1_0(\Om)$. Since $T_m(u_n)$ is bounded in $H^1_0(\Om)$ too, we get $u_n$ bounded in $H^1_0(\Om)$. Taking into account (\ref{codTm}), (\ref{codln}) we have that $u_n$ converges weakly to $u$ in $H^1_0(\Om)$. By the Rellich-Kondrachov's compactness theorem, we can now easily pass to the limit in (\ref{pbhom}) in the distributional sense to deduce that $u=u_0$, the solution of (\ref{pbli}).\par\medskip\noindent
{\it Step 3.} In this and the following step we assume $N\geq 3$, $|E_n|^N$ not necessarily equi-integrable.\par 
For $\delta>0$, we define $w_{n,\delta}$ as the solution of
\beq\label{defwnd} \left\{\ba{l}\dis -{\rm div}\big(A\nabla w_{n,\delta}+E_nw_{n,\delta}\big)+\delta |w_{n,\delta}|^{4\over N-2}w_{n,\delta}=f_n\ \hbox{ in }\Om\\ \ecart\dis w_{n,\delta}\in H^1_0(\Om).\ea\right.\eeq
Using $w_{n,\delta}$ as test function in (\ref{defwnd}), we have
\beq\label{tecoes}\into A\nabla w_{n,\delta}\cdot\nabla w_{n,\delta}\,dx+\into w_{n,\delta} E_n\cdot \nabla w_{n,\delta}\,dx+\delta\into |w_{n,\delta} |^{2N\over N-2}dx=\langle f_n,w_{n,\delta}\rangle.\eeq
\par
In the second term of this equality we use Young's inequality with exponents $2N/(N-2)$, $N$ and 2, to get
$$\ba{l}\dis\left|\into w_{n,\delta} E_n\cdot \nabla w_{n,\delta}\,dx\right|=
\left|\into \Big(\delta^{N-2\over 2N}w_{n,\delta}\Big) \Big({2^{1\over 2}\over \alpha^{1\over 2}\delta^{N-2\over 2N}}E_n\Big)\cdot \Big({\alpha^{1\over 2}\over 2^{1\over 2}}\nabla w_{n,\delta}\Big)dx\right|\\ \ecart\dis
\leq {N-2\over 2N}\delta\into|w_{n,\delta}|^{2N\over N-2}dx+{2^{N\over 2}\over N\alpha^{N\over 2}\delta^{N-2\over 2}}\into|E_n|^Ndx+{\alpha\over 4}\into|\nabla w_{n,\delta}|^2dx.
\ea$$\par
In the last term, Young's inequality also gives
$$\big|\langle f_n,w_{n,\delta}\rangle\big|\leq \|f_n\|_{H^{-1}(\Om)}\|w_{n,\delta}\|_{H^{-1}(\Om)}\leq  {1\over \alpha}\|f_n\|^2_{H^{-1}(\Om)}+{\alpha\over 4}\into |\nabla w_{n,\delta}|^2dx.$$
Using also (\ref{elip}) in the first term, we deduce from (\ref{defwnd})
\beq\label{estwnd}{\alpha\over 2}\into\big|\nabla w_{n,\delta}\big|^2dx+{N+2\over 2N}\delta  \into |w_{n,\delta}|^{2N\over N-2}dx\leq {2^{N\over 2}\over N\alpha^{N\over 2}\delta^{N-2\over 2}}\into|E_n|^Ndx+{1\over \alpha}\|f_n\|^2_{H^{-1}(\Om)}. \eeq
Thus, for every $\delta>0$, $w_{n,\delta}$ is bounded in $H^1_0(\Om)$. Thanks to Rellich-Kondrakov's compactness theorem this allows us to pass to the limit in (\ref{defwnd})
 in the distributional sense,  to deduce
\beq\label{cownd} w_{n,\delta}\rightharpoonup w_{\ast,\delta}\ \hbox{ in }H^1_0(\Om),\eeq
with $w_{\ast,\delta}$ the solution of
\beq\label{defwd} \left\{\ba{l}\dis -{\rm div}(A\nabla w_{\ast,\delta}+E_0w_{\ast,\delta})+\delta |w_{\ast,\delta}|^{4\over N-2}w_{\ast,\delta}=f\ \hbox{ in }\Om\\ \ecart\dis w_{\ast,\delta}\in H^1_0(\Om).\ea\right.\eeq
\par
We take $u_0$ the solution of (\ref{pbli}). Taking $T_\rho(w_{\ast,\delta}-u_0)$, with $\rho>0$ as test function in the difference of (\ref{defwd}) and (\ref{pbli}), we have
$$\ba{l}\dis \int_{\{|w_{\ast,\delta}-u_0|<\rho\}}\hskip-10pt \big(A\nabla (w_{\ast,\delta}-u_0)+ E(w_{\ast,\delta}-u_0)\big)\cdot\nabla (w_{\ast,\delta}-u_0)\,dx\\ \ecart\dis+\delta\into\big( |w_{\ast,\delta}|^{4\over N-2}w_{\ast,\delta}-|u_0|^{4\over N-2}u_0\big)T_\rho(w_{\ast,\delta}-u_0)\,dx=-\delta \into |u_0|^{4\over N-2}u_0T_\rho(w_{\ast,\delta}-u_0)\,dx.\ea$$
Using Young's inequality in the second term of the first integral and
$$\exists c_N>0:\quad c_N\big||x|+|y|\big|^{4\over N-2}|x-y|^2\leq \big(|x|^{4\over N-2}x-|y|^{4\over N-2}y\big)(x-y),\quad\forall\,x,y\in\RR,$$
in the second integral, we get
\beq\label{estwd0}\ba{l}\dis{\alpha\over 2}\int_{\{|w_{\ast,\delta}-u_0|<\rho\}}\hskip-16pt |\nabla (w_{\ast,\delta}-u_0)|^2dx+c_N\delta \into \big(|w_{\ast,\delta}|+|u_0|\big)^{4\over N-2}|w_{\ast,\delta}-u_0||T_\rho(w_{\ast,\delta}-u_0)|\,dx\\ \ecart\dis\leq {1\over 2\alpha}\int_{\{|w_{\ast,\delta}-u_0|<\rho\}}\hskip-16pt |E|^2|w_{\ast,\delta}-u_0|^2dx+
\delta \into |u_0|^{N+2\over N-2}|T_\rho(w_{\ast,\delta}-u_0)|\,dx.\ea\eeq
This proves that  for every $\rho>0$, $T_\rho(w_{\ast,\delta}-u_0)$ is bounded in $H^1_0(\Om)$. Dividing by $\rho$ and taking the limit when $\rho$ tends to zero, thanks to the Lebesgue dominated convergence theorem, we get
$$\into \big(|w_{\ast,\delta}|+|u_0|\big)^{4\over N-2}|w_{\ast,\delta}-u_0|dx\leq {1\over c_N}\into |u_0|^{N+2\over N-2}dx,$$
and thus $w_{\ast,\delta}$ is bounded in $L^{N+2\over N-2}(\Om)$.
Hence, for a subsequence of $\delta$ which converges to zero, still denoted by $\delta$, there exits $w_\ast$ such that 
\beq\label{conwel} \left\{\ba{l}\dis w_{\ast,\delta}\to w_\ast\ \hbox{ a.e. in }\Om,\ \ w_{\ast,\delta}\rightharpoonup w_\ast\ \hbox{ in }L^{N+2\over N-2}(\Om)\\ \ecart\dis
T_\rho(w_{\ast,\delta}-u_0)\rightharpoonup T_\rho(w_\ast-u_0)\hbox{ in }H^1_0(\Om),\ \forall\,\rho>0.\ea\right.\eeq
By the lower semicontinuity of the norm in $H^1_0(\Om)$, this allows us to pass to the limit when $\delta$ tends to zero in (\ref{estwd0}) to deduce 
$$\int_{\{|w_\ast-u_0|<\rho\}}\hskip-10pt |\nabla (w_\ast-u_0)|^2dx\leq {1\over \alpha^2}\int_{\{|w_\ast-u_0|<\rho\}}\hskip-10pt |E|^2 |w_\ast-u_0|^2dx.$$
Dividing by $\rho^2$ and taking the limit when $\rho$ tends to zero, this proves
$${1\over \rho}T_\rho(w_\ast-u_0)\rightarrow 0\ \hbox{ in }H^1_0(\Om)\ \hbox{ when }\rho\to 0.$$
Combined with 
$${1\over \rho}T_\rho(w_\ast-u_0)\rightharpoonup {\rm sgn}(w_\ast-u_0)\ \hbox{ a.e. in }\Om,$$
we get
\beq\label{weqz} w_\ast=u_0\ \hbox{ a.e. in }\Om.\eeq\par
\par\medskip\noindent
{\it Step 4.} For $\rho>0$  we take, similarly to the Step 3, $T_\rho(w_{n,\delta}-u_n)$ as test function in the difference of  (\ref{defwnd}) and (\ref{pbhom}). We get
$$\ba{l}\dis{\alpha\over 2}\int_{\{|w_{n,\delta}-u_n|<\rho\}}\hskip-12pt |\nabla (w_{n,\delta}-u_n)|^2dx+\delta \into |w_{n,\delta}|^{4\over N-2}w_{n,\delta} T_\rho(w_{n,\delta}-u_n)\,dx\\ \ecart\dis\leq {1\over 2\alpha}\int_{\{|w_{n,\delta}-u_n|<\rho\}}\hskip-12pt |E_n|^2|w_{n,\delta}-u_n|^2dx.\ea$$
Taking into account (\ref{codTm0}) and (\ref{cownd}), and defining ${\cal E}$ by (it exists for a subsequence)
$$|E_n|^2\rightharpoonup {\cal E}\ \hbox{ in }L^{N\over 2}(\Om),$$
we can  pass to the limit in $n$ in this inequality by semicontinuity and the Rellich-Kondrachov's compactnes theorem to deduce
$${\alpha\over 2}\int_{\{|w_{\ast,\delta}-u|<\rho\}}\hskip-12pt |\nabla (w_{\ast,\delta}-u)|^2dx+\delta \into |w_{\ast,\delta}|^{4\over N-2}w_{\ast,\delta} T_\rho(w_{\ast,\delta}-u)\,dx\leq {1\over 2\alpha}\int_{\{|w_{\ast,\delta}-u|\leq \rho\}}\hskip-12pt {\cal E}|w_{\ast,\delta}-u|^2dx.$$\par
Now we pass  to the limit when $\delta$ tends to zero thanks to (\ref{conwel}) and (\ref{weqz}), to get
$$\int_{\{|u_0-u|<\rho\}}\hskip-10pt |\nabla (u_0-u)|^2dx\leq {1\over \alpha^2}\int_{\{|u_0-u|\leq\rho\}}\hskip-10pt {\cal E}|u_0-u|^2dx$$
Dividing by $\rho^2$ and passing to the limit when $\rho$ tends to zero we deduce  as at the end of Step 3  that $u=u_0$. This finishes the proof. \cqfd\par\bigskip
\begin{Rem} One of the applications of Theorem \ref{Thgho} is the existence of solutions for some control problems in the coefficients. In this way, combined with Fatou's Lemma it immediately proves the existence of solution for
$$\ba{c}\dis \min\left\{\into \big(G(x,u)+\mu|E|^p\big)dx:\ E\in L^p(\Om)^N\right\}\\ \ecart\dis
\left\{\ba{l} -{\rm div}(A\nabla u+Eu)=f\ \hbox{ in }\Om\\ \ecart\dis u\in H^1_0(\Om),\ea\right.\ea$$
with $p>2$ if $N=2$, $p\geq N$ if $N>2$, $f\in H^{-1}(\Om)$ and $G:\Om\times \RR\to\RR$  such that
$$G(.,s)\ \hbox{ measurable, }\forall\, s\in\RR,\quad G(x,.)\hbox{ continuous for a.e. }x\in\Om,$$
$$\exists a\in\RR, b\geq 0, \ \hbox{ such that }\ G(x,s)\geq a-b|s|,\quad\forall\,s\in\RR,\ \hbox{ a.e. }x\in\Om.$$
\end{Rem}
\section*{Aknowledgments} This work has been partially supported by the project PID2020-116809GB-I00
of the {\it Ministerio de Ciencia e Innovaci\'on}  of the government of Spain.

\end{document}